\documentclass[12pt]{amsart}

\usepackage{amsgen,amsthm,amstext,amsbsy,amsopn,amsfonts,amssymb,enumerate}

\usepackage{soul}
\usepackage{hyperref}
\usepackage{amsmath}
\usepackage[usenames]{xcolor}

\newcommand\la{\langle}
\newcommand\ra{\rangle}

\newcommand\sg{{\mathfrak{s}}}

\newcommand\HH{\mathbb H}

\newcommand\RR{\mathbb R}

\newcommand\ad{\operatorname{ad}}

\newcommand\GL{\operatorname{GL}}

\theoremstyle{plain}

\newtheorem{thm}{Theorem}[section]

\theoremstyle{definition}

\newtheorem{rem}[thm]{Remark}
\newtheorem{example}[thm]{Example}

\begin{document}

\title[The mean curvature flow on solvmanifolds]{The mean curvature flow on solvmanifolds}

\author[Arroyo]{Romina M. Arroyo}
\address{FAMAF \& CIEM - Universidad Nacional de C\'ordoba. Av  Medina Allende s/n - Ciudad Universitaria.
ZIP:X5000HUA - Cordoba, Argentina}
\email{arroyo@famaf.unc.edu.ar}

\author[Ovando]{Gabriela P. Ovando}
\address{Departamento de Matem\'atica, ECEN - FCEIA, Universidad Nacional de Rosario. Pellegrini 250, 2000 Rosario, Santa Fe, Argentina}
\email{gabriela@fceia.unr.edu.ar}

\author[Perales]{Raquel Perales}
\address{Institute of Mathematics at the National Autonomous University of Mexico. Leon 2, altos. Oaxaca, Oaxaca.  68000. Mexico}
\email{raquel.perales@im.unam.mx}

\author[S\'aez]{Mariel S\'aez}
\address{Facultad de Matem\'atica, P. Universidad Cat\'olica de Chile, Av. Vicuna Mackenna 4860, 690444 Santiago, Chile.}
\email{mariel@mat.puc.cl}

\thanks{{\it (2000) Mathematics Subject Classification}:  53A10, 53C42, 53C44, 22E25, 22F30}




\begin{abstract}  
This work is a survey of  the most relevant background material to motivate and understand the construction and classification of translating solutions to mean curvature flow on a family of solvmanifolds. 

We introduce the mean curvature flow and some known results in the field. In particular we explore the notion of translating solution in the Euclidean space and extensions into other Riemannian manifolds. We also include a discussion on solvmanifolds and some elements of its geometry that are relevant to our work. 
 We finish by posing the equations that describe translating solutions to mean curvature flow on our family of 3-dimensional solvmanifolds with some additional assumptions.
 
This project emerged at the ``Latin American and Caribbean Workshop on Mathematics and Gender'' held at Casa Matem\'atica Oaxaca (CMO) from May 15-20, 2022.

\end{abstract}

\maketitle

\setcounter{tocdepth}{1}
\tableofcontents

\section{Introduction}\label{sec:intro}

The mean curvature flow is an extrinsic geometric flow with an extensive history that includes applied problems as well as purely geometric questions. This equation describes the gradient flow of the area functional, describing the deformation of submanifolds as they minimize their volume in the steepest direction. The mean curvature flow was first described in the 1950's  by Mullins \cite{Mullins} in the context of material science, where the motion of interfaces between phases of a material can be described by this equation due to the minimization of certain energy. The first mathematically rigorous work on the mean curvature flow was carried out by Brakke in 1978 \cite{Brakke}, where he used the tools of Geometric Measure Theory to  define solutions to this equation. In the 1980's, Huisken  developed fundamental aspects of the theory for classical solutions \cite{Huisken, Huisken1986ContractingCH}. Subsequently, there have been many results regarding the mean curvature flow (see \cite{ABGL} and references therein).\\

Solutions to mean curvature flow with additional structure have played an important role in the theory, providing explicit examples and serving as models for the behavior of more general solutions to the equation. To compute explicit solutions in $\mathbb{R}^n$, symmetries of Euclidean space, such as rotational symmetry, translation invariance, or graphicality, are usually exploited.  A particularly important subset of these specific solutions is the category of solitons or self-similar solutions, which we will introduce in Section \ref{MCF} and elaborate on their relevance. \\

While the mean curvature flow is well-posed on any Riemannian manifold, there are not many works adapted to specific non-Euclidean geometries. Explicit examples of solutions must be considered in each special case. Nevertheless, several works suggest that the geometry of submanifolds in non-Euclidean spaces can exhibit richness and diversity compared to its Euclidean counterpart. An illustrative case study is the minimal surface theory, which corresponds to the static case of mean curvature flow, in $\mathbb{H}^2\times\mathbb{R}$ (for a non exhaustive list of references, see  \cite{Daniel, FMMR, Hau, NR}). Notably, the manifold $\mathbb{H}^2\times\mathbb{R}$ belongs to a broader family of 3-manifolds known as solvmanifolds. These can be described as $\mathbb{R}^3$ with a family of metrics given by a warped product and a group structure, playing a crucial role in Thurston's geometries and the topological classification of 3-manifolds. \\

The role of solvmanifolds in classifying 3-manifolds establishes a connection between them and general 3-manifolds through Ricci flow. This process breaks down general manifolds into fundamental pieces. On the other hand, several authors [17, 37, 48] have studied the interplay between mean curvature flow and Ricci flow.
It would be certainly interesting to understand whether these results  may provide a tool to connect results of mean curvature flow in solvmanifolds with more general 3-dimensional geometries. \\

In this work we are interested in studying solitons to mean curvature flow on 3-dimensional solvmanifolds. An advantage of working in this setting, as opposed to a more general Riemmanian manifold, is that symmetries of the ambient space can be explicitly stated.  Specific examples exploiting these symmetries have been already explored by some authors in some special situations, and we will discussed this in Section \ref{knownexamples}.\\

The main goal of this survey is to provide sufficient background to motivate and pose the question of constructing  {\it translating solutions} to mean curvature flow on general 3-dimensional  solvmanifolds. To this end, the structure of this manuscript  is as follows: In Section \ref{MCF}, we motivate and set up the mean curvature flow on Riemannian manifolds. We also explain the concept of {\it translating solutions} in Euclidean space and how this is extended to other geometries. In Section \ref{sec:LieGroup}, we define a particular family of  3-dimensional  solvmanifolds and describe their Killing vector fields. In Section \ref{knownexamples} we discuss several known results. Finally,  in Section \ref{extensions}  we propose a setup to extend some results from the previous section.  While we expect that readers with a basic background in Riemannian geometry can follow our computations, we include details for their convenience.

\section{The Mean Curvature Flow}\label{MCF}
The (co-dimension 1) mean curvature flow is a geometric equation that considers the evolution in time  of an  embedding $\varphi: M^n\times[0,T)\to \mathcal{N}^{n+1}$  of an $n$--manifold $M$ into an $(n+1)$--dimensional Riemannian manifold $(\mathcal{N},g)$, with normal speed equal to its mean curvature. This translates into 
 the evolution equation
 \begin{equation} 
 g\left(\frac{\partial \varphi}{\partial t}, \nu\right)=-H, \label{eq:MCF}
  \end{equation}
where

$H$ and $\nu$  are the (scalar) mean curvature and 
a consistent choice of unit normal vector field of $M$, respectively. \\

\subsection{ The Euclidean case} When  $\mathcal{N}=\mathbb{R}^{n+1}$  with the standard metric is considered
this equation  has been extensively studied (see the introduction above for references). This situation is referred to as the ``classical case'' in this manuscript.  Here, we list a few well-known examples of explicit solutions in $\mathbb{R}^{n+1}$.

\medskip

\begin{itemize}
\item Spheres of radius $R(t)=\sqrt{R_0-2nt}$ and cylinders of radius $R(t)=\sqrt{R_0-2mt}$ with $m\leq n$ and $n-m$ flat directions.
\item For $n=1$, the grim reaper $\varphi(x,t)=(x, t-\ln \cos x)$ and for general $n$ we have $\varphi(x,t)\times \mathbb{R}^{n-1}$.
\item The bowl soliton: This is a rotationally symmetric graph that exists for dimensions bigger or equal than 2, and it is invariant by translation in time.
\end{itemize}

Note that the first example is only defined until the finite time $T=\frac{R_0}{2nt}$, while the other two examples are defined for all times. \\

A PDE argument that uses the maximum principle  allows us to show that evolving surfaces that are disjoint at the initial time remain disjoint as long as the evolutions are smooth (see Theorem 6.16 in \cite{ABGL}). This, in particular, implies that all compact initial conditions develop singularities in finite time (since they can be enclosed by a large sphere which becomes extinct in finite time), and the understanding of these singularities is one of the main goals in the field.\\

It is not difficult to show (using arguments from PDE theory, see Theorem 6.20 in \cite{ABGL}  for instance) that singularities occur when the second fundamental form of the evolving surface is unbounded. One possible classification of singularities is according to the ``rate'' at which this blow up occurs. More precisely, if the norm of the second fundamental form is of order $\frac{1}{\sqrt{T-t}}$
(where $T$ is the maximal time of existence), then the singularity is labeled as a ``type I singularity''. Otherwise, singularities are known as 
``type II".\\
 
To gain a better understanding of  singularities, rescaling procedures adapted to each type of singularity,  often referred to  as ``blow-ups'', have been proposed (see \cite{Mante} for a precise description).  In the limit these blow ups converge to a special solution of the mean curvature flow.
In the case of type I singularities, the limit of such rescaling converges to a {\em shrinking self-similar solution}, which is a solution of the form $\lambda(t) M_0$, where $M_0$ is a fixed manifold and $\lambda '(t)<0$. Under additional convexity assumptions, these solutions are seen to be spheres and cylinders. For a classification of shrinking self-similar solutions under certain assumptions, see \cite{AL, Huisken1990AsymptoticbehaviorFS, Huisken1993LocalAG}.\\

On the other hand, for type II singularities, Hamilton proposed a rescaling procedure that ensures that the limit is an {\it eternal solution} -- a solution that exists for all times (see \cite{Mante} for a precise description). Translating self-similar solutions (or translators) are a particular case of eternal solutions that can be expressed as $M_t=M_0+tv$, where $v$ is the direction of translation. 
It is a standing question whether it is always possible to perform rescalings of type II singularities that produce translators. Under some additional assumptions, this is the case (see the work of Choi, Haslhofer and Hershkovits in \cite{ChoiHaslhoferJershkovits} for example), 
motivating the classification of translators
 (see for instance \cite{AL, BLT, HIMW, HMW}.)\\

Another topic of interest in the field has been the study of entire graphs evolving by mean curvature flow.  In contrast with the compact case, these solutions do not develop singularities in finite time. The seminal work of Ecker and Huisken in \cite{EH1, EH2} shows long-time existence for any locally Lipschitz condition, along with several relevant estimates and other results. In this case,  {\em expanding self-similar solutions} (or expanders) play a role. These solutions are described as  $\lambda(t) M_0$, where $M_0$ is a fixed manifold and $\lambda'(t)>0$. \\

\subsection{Self-similar solutions in general Riemannian manifolds}
As can be seen from the previous discussion,  the concept of self-similar solution (shrinking, translating or expanding) depends on the symmetries of Euclidean space. An extension considered in the literature of this concept  to other manifolds  involves solutions $\varphi_t(M_0)$, where $\varphi$ is the flow associated to a vector field $V$ that is either conformal or Killing (see for example \cite{CoMaRi, KO,LO,Pi1, Pi2}). Recall that $V$ is called conformal if it satisfies $\mathcal{L}_Vg= 2f g$, where $\mathcal{L}$ is the Lie derivative, $g$ the ambient metric and $f$ is called the potential function. In the special case, when $f\equiv 0$ the vector field $V$ is Killing, which in the Euclidean setting is identified with translators. On the other hand, non trivial conformal fields can be associated (in the Euclidenan setting) to homothetic self-similar solutions (that is, shrinkers and expanders).

\section{The solvmanifolds $(S,g)$}\label{sec:LieGroup}

The aim of this section is to introduce the family of $3$-dimensional solvmanifolds for which we will study the Mean Curvature Flow. To achieve this goal and facilitate self-contained reading, we will briefly review the notion of a Lie group equipped with a left-invariant metric.

\subsection{Lie groups, solvmanifolds and groups of isometry.}
A Lie group  is a group $G$ endowed  with a differentiable structure, for which  the product map $G\times G \to G$ and the inverse map $G\to G$ are differentiable maps. 
So,  left-translations by elements of the group,  $L_p(g)
=pg$  for every $g\in G$, are diffeomorphisms of the Lie group. \\

A vector field $X\in \mathcal X (G)$ is called left-invariant if $dL_p \circ X=X \circ L_p$ for every $p\in G$. The set of left-invariant vector fields of $G$ forms a finite dimensional subspace closed under the Lie bracket.  Indeed, every left-invariant vector field of $G$ is determined by its value at the identity element. This subspace is called the {\em Lie algebra} of $G$ and will be denoted by $\mathfrak g$.\\

Conjugations maps $I_g:G\to G$ are automorphisms of $G$ defined  by
$$I_g(x)=gxg^{-1},$$
and are diffeomorphisms of $G$. Since $I_g(e)=e$ for all $g\in G$, it follows that the differential 
$d(I_g)_e:T_eG\to T_eG$ induces an isomorphism at the Lie algebra level. It is not hard to verify that $Ad:G\to \GL(\mathfrak g)$, given by  $Ad(g)=d(I_g)_e$, is a homomorphism of Lie groups, that is, a representation, which is called the adjoint representation.\\

A Lie group $G$ is called {\em unimodular} if and only if 
$$\det Ad(g) = 1, \qquad\mbox{ for all } \quad  g \in G.$$
 At the Lie algebra level this implies that  $\mathrm{tr} \,  ad(X) =0$ for any $X \in \mathfrak g$,  with equivalence whenever the Lie group $G$ is connected. Here $ad(X):\mathfrak g \to \mathfrak g$ denotes the map $ad(X)(Y):=[X,Y]$ and $\mathrm{tr}(\cdot)$ denotes the trace operator.\\

At the Lie algebra $\mathfrak g$ one has the {\em derived series} defined by:
$$D^0(\mathfrak g)=\mathfrak g, \quad  D^j(\mathfrak g)=[D^{j-1}(\mathfrak g), D^{j-1}(\mathfrak g)] \quad\mbox{ for }j\geq 1.$$
If there is $k>0$ such that $D^k(\mathfrak g)=0$, the Lie algebra is called {\em solvable} and the corresponding Lie group is called a solvable Lie group. \\

One says that a  group $G$ acts (on the left) on a set $M$ if there is a function $\mu:G\times M\to M$ such that:
\begin{itemize}
	\item $\mu(gh,x)=\mu(g, \mu(h,x))$ for all $g,h\in G, x\in M$. 
	\item $\mu(e,x)=x$ for every $x\in M$, where $e$ denotes the identity element in $G$.
\end{itemize}
Given an action,  the {\em orbit} of $x\in M$ is given by:
$$G\cdot x=\{\mu(g,x)\mbox{ for }g\in G\},$$
and the isotropy subgroup at $x\in M$ is given by:  
$$G_x=\{g\in G\,:\, \mu(g,x)=x\}.$$
We note that if the action of $G$ on $M$ is transitive all the isotropy groups are conjugated and isomorphic.

\begin{example}
Let $G$ be a Lie group. 
\begin{itemize}
\item $G$ acts on itself by translations on the left: $G\times G \to G$ given by $(p,x)\to px$. This action has trivial isotropy.
\item Another action of $G$ on $G$ is by conjugation $G\times G\to G$, $(p,x)\to pxp^{-1}.$ In this case the isotropy subgroup at $x\in G$,
$G_x$,  is 
	usually called the centralizer of $x$ in $G$. 
\end{itemize}	
	\end{example}

If a differentiable manifold $M$ admits a transitive action of a Lie group $G$, one says that the manifold is a homogeneous space (more details on homogeneous spaces can be read in \cite{Hel}). A {\em solvmanifold} is a connected Riemannian
manifold which is acted upon transitively by a solvable Lie group of isometries (see \cite{GW88}).

\subsubsection{Isometry groups of  $3$-dimensional Lie groups $(G,g)$}

Let $G$ be a  Lie group with  a left-invariant metric $g$ defined on it. Every isometry  can be written as $f=L_p \circ \varphi$, a composition of an isometry $\varphi$ that preserves the identity element of $G$ and a translation on the left $L_p$,  by some element $p\in G$.  More explicitly,  if $f$ is an isometry of $G$ we can assume $f(e)=p$ for some $p \in G$. Then, $L_{p^{-1}}\circ f$ is an isometry of $G$ preserving the identity element. Thus, $L_{p^{-1}}\circ f=\varphi$. Equivalently,  $f=L_p\circ \varphi$ with $\varphi$ an isometry preserving the identity element. \\

The subgroup of translations on the left is isomorphic to the Lie group. In fact $L:G \to \mathrm{Isom}(G,g)$ given by $g\to L_g$ is a  monomorphism onto its image.\\

Thus, the isometry group of $(G,g)$, denoted by $\mathrm{Isom}(G,g)$, is a group with the composition of functions as the group operation and admits a decomposition as a product of subgroups: 
$$\mathrm{Isom}(G,g)= G \,.\, \mathrm{Isom}(G,g)_e.$$

This is a generalization of the situation in Euclidean space, where every isometry decomposes as $T_a\circ R$, with $T_a$ representing a translation by $a$ (defined as $T_a(x)=x+a$) and $R$ being an orthogonal map. Notice that the sum of vectors in Euclidean space is commutative, a property that does not necessarily hold in Lie groups. \\

The isometry groups of simply connected Lie groups of dimension three with a left-invariant metric were computed in \cite{HL} and \cite{CR} for the unimodular and non-unimodular cases, respectively.\\

There exist, up to Lie group isomorphisms, three simply connected unimodular solvable Lie groups of dimension three: 
 $Nil_3$, $Sol_3$ and  $\widetilde{E}_0(2)$.  The isometry groups of these were obtained  in 
 Theorem 3.2 for $Nil_3$, Theorem 3.3 for $Sol_3$ and  Theorem 3.4 for $\widetilde{E}_0(2)$ in  \cite{HL}.
\\ 

For the non-unimodular setting, every three-dimensional solvable Lie algebra decomposes as a sum of vector spaces $\mathbb R^2\oplus \mathbb R$, where $\mathbb R^2$ is an abelian ideal spanned by vectors $V_1, V_2$ and the action of $V_3$ on $\mathbb R^2$ is determined by a $2\times 2$ real matrix $A$, with $\mathrm{tr} (A)=2$, whose determinant $c=\det(A)$ is an isomorphism invariant, except when $A$ is the identity. This invariant was introduced by Milnor in \cite{Mi}. The action of $V_3$ is given either by
\begin{itemize}
	\item the identity or
	\item a matrix of the form $\left(\begin{matrix}
		0 & 1 \\-c & 2
	\end{matrix} \right)$.
\end{itemize}

For $A=\mathrm{Id}$ the isometry group was computed in \cite[Theorem 3.1]{CR}.\\

For every $c$ as above, the corresponding Lie group is denoted by $G_c$. The isometry groups were computed  
\begin{itemize}
\item in Theorem 3.2 in \cite{CR} for $G_0$, 
\item in Theorem 3.5 in \cite{CR} for $G_1$  and 
\item in Theorems 3.6, 3.7 and 3.8 \cite{CR} for $G_c$, with $c \neq 0, 1$.
\end{itemize}

\subsection{The solvmanifold $(S,g)$}

Three dimensional solvmanifolds have been  extensively used for various purposes.  In particular,  they give models for three-dimensional compact spaces in the  geometrization conjecture. We refer the reader to the foundational work of John Milnor \cite{Mi} for  properties of these spaces (such as curvature properties).
We now describe the family of solvmanifolds we are interested in.\\

The three-dimensional solvmanifolds we are interested in are constructed as follows. 
Consider the underlying differentiable structure of $\RR^3$ with its usual coordinates and the group structure defined by:
$$(x_1,y_1,z_1)(x_2,y_2,z_2)=(x_1+e^{\lambda_1 z_1}x_2, y_1+e^{\lambda_2 z_1}y_2, z_1+z_2),$$
where $\lambda_1, \lambda_2\in \RR$. 
Note that for $\lambda_1=\lambda_2=0$, $S$ is isomorphic as a Lie group to $\RR^3$.  
To avoid this case, which has been extensively studied, we will always assume that $\lambda_1 \neq 0$.\\

In this product, the identity element is $e=(0,0,0)$, and for every  $(x,y,z)\in \RR^3$, its inverse is given by $(x,y,z)^{-1}=(-e^{-\lambda_1 z}x, -e^{-\lambda_2z}y, -z)$. \\

A basis of left-invariant vector fields on $S$ is obtained as follows:
\begin{equation}\label{leftinvvectorfields}
\begin{split}
E_1(x,y,z) :=\frac{d}{ds}|_{s=0} (x,y,z)(s,0,0)=e^{\lambda_1z}\partial_x,\\
E_2(x,y,z) :=\frac{d}{ds}|_{s=0} (x,y,z)(0,s,0)=e^{\lambda_2z}\partial_y,\\
E_3(x,y,z):=\frac{d}{ds}|_{s=0} (x,y,z)(0,0,s)=\partial_z.
\end{split}
\end{equation}
This basis satisfies the following non-trivial Lie bracket relations:
$$[E_3,E_1]=\lambda_1 E_1, \qquad [E_3,E_2]=\lambda_2 E_2.$$
From these bracket relations, we see that $\sg$ is solvable. In fact, its derived series ends in the second step: $\sg \supseteq [\sg, \sg] \supseteq [[\sg, \sg], [\sg, \sg]]=0$. \\

Take the Riemannian metric on $S$ for which the basis mentioned above is orthonormal; that is, at each $(x,y,z)$:
\begin{equation}
g=e^{-2\lambda_1z}dx^2 + e^{-2\lambda_2z}dy^2 + dz^2. \label{sol metric} 
\end{equation}
For $g$ as above, left-translations $L_p$ are isometries, $p\in S$. Since this action of $S$ by isometries is transitive, the manifold $S$ is a homogeneous space.  In particular, this solvmanifold is complete (Theorem 4.5 Ch. IV in \cite{KN}).\\

Let $\nabla$ denote the Levi-Civita connection on $S$. The Koszul formula for left-invariant vector fields $U,V$ gives
$$2 \nabla _U V = [U,V]- \ad^t(U)(V) - \ad^t(V)(U),  \
$$
where $\ad^t(U)$ denotes the transpose of $\ad(U):\sg\to \sg$ with respect to the metric $g$. From it, 
one gets

\begin{equation}\label{connection}
\begin{split}
\nabla_{E_1}E_1= \lambda_1 E_3,  \quad   \nabla_{E_1}E_3= -\lambda_1 E_1, \\
\nabla_{E_2}E_2= \lambda_2 E_3,  \quad   \nabla_{E_2}E_3= -\lambda_2 E_2,
\end{split}
\end{equation}
and zero otherwise.\\

\begin{rem}

For $\lambda_1=\lambda_2$, it is not hard to compute that the sectional curvatures of $(S,g)$ are equal to
$$sec(v,w)=-\lambda_1^2,$$
for orthonormal vectors $v,w$ in $T_pS$.
Hence, $(S, g)$ as a simply connected manifold of constant sectional curvature is isometric to the hyperbolic space (see, for instance, Theorem 11.12, Ch. 11 in \cite{Lee}).\\

Furthermore,  there exists an isomorphism of Lie groups between $S$ and the $3$-dimensional hyperbolic space.
To see this, we recall the following explicit construction: \\

The group structure on a half-space model of the real hyperbolic space 
\[
\mathbb H^n= \{  (x_1,\hdots , x_{n-1}, t) \in  \Bbb R^n \, | \, t > 0\}
\]
is given by
$$(x_1, \hdots , x_{n-1},t) (y_1, \hdots , y_{n-1}, s) = \bigl( x_1 + ty_1, \hdots  , x_{n-1} + ty_{n-1}, ts
\bigr).$$ 
This is a Riemannian manifold with constant negative curvature when endowed with the metric 
\[
g_{\lambda}  = \frac{1}{x_n^2}\bigl(dx_1^2 + \hdots + dx^2_{n-1} + \lambda dx^2_n\bigr)
\] for $\lambda  \neq  0$. Then the map $\psi:S \to \mathbb H^3$ given by
$$ \psi :  (x_1,x_2 , x_3) =  (x_1,x_2 ,e^{x_3})$$
is an isomorphism of Lie groups. 
\end{rem}

\bigskip

We now describe $\mathrm{Isom}(S,g)$. First,  we have that left-translations on $S$ are generated by the following 
one-parameter groups:
$$\begin{array}{c}
(x,y,z) \to (x+c, y, z),  \\
	(x,y,z) \to (x, y+c, z),  \\ 
 (x,y,z) \to (e^{\lambda_1 c}x, e^{\lambda_2c}y, z+c),
	\end{array}
	$$
	where $c\in \mathbb{R}$.


In view of the scalars $\lambda_i$, the unimodular cases correspond to the condition $\lambda_1+\lambda_2=0.$ These are called the $Sol_3$ models. The MCF in the solvmanifold $Sol_3$ (with $\lambda_1=1$) was studied in \cite{Pi1}, where the corresponding isometry group was considered.  \\

We focus here in the non-unimodular situation, when equipped with the left-invariant metric in Equation \eqref{sol metric}.
In this case, from the results in \cite{CR} mentioned above, the isometry groups are
\begin{itemize}
	\item For $\lambda_1=\lambda_2$, the corresponding isometry group is\\ $\mathrm{Isom}(S,g)=O(3,1)=S . O(3)$.
	\item For $\lambda_2=0$ we have $\mathrm{Isom}(S,g)=S . SO(2)$.
		\item For the other cases, $\mathrm{Isom}(S,g)=S$, that is, it consists of left-translations. 
	\end{itemize}

It is worth noting that the following map is an orientation reversing isometry fixing the identity element on $S$:
$$(x,y,z) \to (-x,y,z).$$ \\

Since we are interested in studying translators, here we list the Killing vector fields of the family $(S,g)$.
Let us note that the subalgebra of Killing vector fields is in correspondence with the Lie algebra of $\mathrm{Isom}(S,g)$.\\

The Killing vector fields corresponding to translations on the left, are the right-invariant vector fields:
\begin{align*}
\widetilde{E}_1(x,y,z) & = \frac{d}{ds}|_{s=0}(s,0,0)(x,y,z)=\partial_x=e^{-\lambda_1 z}E_1,\\
\widetilde{E}_2(x,y,z) & =  \partial_y= e^{-\lambda_2 z}E_2,\\
\widetilde{E}_3(x,y,z) & =  \lambda_1 x \, \partial_x + \lambda_2 y \, \partial_y+ \partial_z= \lambda_1 x e^{-\lambda_1 z}E_1 + \lambda_2 y e^{-\lambda_2 z}E_2+ E_3.
\end{align*}

\bigskip

In the cases of isometries fixing the identity element of $S$ we also have the following:

\begin{itemize} 

  \item For the case $\lambda_1=\lambda_2$, remember that a vector field $\widetilde{E}$ is Killing if and only if $\la \nabla_Y \widetilde{E},Z \ra = - \la Y, \nabla_Z \widetilde{E} \ra $, for all $Y,Z \in \mathcal X (S)$. Then, by using (\ref{leftinvvectorfields}) and (\ref{connection}), it is easy to see that the following vector fields are also Killing:
\begin{align*}
\widetilde{E}_4(x,y,z) & = - y e^{-\lambda_1 z}E_1 + x e^{-\lambda_1 z}E_2,\\
\widetilde{E}_5(x,y,z) & =  y E_2 - \frac{1}{2} y^2 \lambda_1 e^{-\lambda_1 z}E_3,\\
\widetilde{E}_6(x,y,z) & =  x E_1 - \frac{1}{2} x^2 \lambda_1 e^{-\lambda_1 z}E_3.
\end{align*}
   
 \item For the case  $\lambda_2=0$ {i.e. $S$ is isometric to the product space $\mathbb H^2 \times \RR$,} one has the  Killing vector field $$\widetilde{E}(x,y,z)=-a z \partial_y + a y \partial_z=-a zE_2+ayE_3,$$ 
for any $a \in \RR$. This corresponds to the following monoparametric group of isometries of  $S$:
$$E_u(x,y,z)=\left(\begin{matrix}
1 & 0 & 0\\
0 & 	\cos(au) & -\sin(au) \\
0 & 	\sin(au) & \cos(au)
\end{matrix} \right)\left(\begin{matrix}
x\\y\\z \end{matrix} \right)  \mbox{ for }u\in\RR,a\in \RR^*.$$
\end{itemize}

%
%

\section{Solitons of Mean Curvature Flow on solvamnifolds} \label{knownexamples}

In this section, we review some known constructions on solvmanifolds. These constructions have been performed for specific choices of the parameters $\lambda_1$ and $\lambda_2$. We subdivide the discussion into each of these cases. We would like to remark that we are interested on constructions related to solitons of mean curvature flow. However, in the literature there are also related results on other flows \cite{LP}, as well as related problems, such as prescribed curvature functions \cite{BO}.

\subsection{Euclidean space: $\lambda_1=\lambda_2=0$}

As previously mentioned, the literature in this case is extensive (see for instance \cite{BLT, HIMW} for a discussion). We will not include all the known results for this case, but  describe the most relevant ones for our discussion.\\

The best known-example is the grim reaper on the plane, which was already described in Section \ref{MCF} and it is the unique translating solution to curve shortening flow on the plane \cite{Giga}. This is a graphical translating solution given by $\varphi(x,t)=(x, t-\ln \cos x)$. For general dimensions,  it can be extended trivially in other directions, i.e. $\varphi(\cdot,t)\times \mathbb{R}^{n-1}$. The grim reaper also gives rise to another solution in higher dimensions,  obtained through a suitable rotation with respect to a two dimensional plane generated by the translating direction and an orthogonal direction. This solution is known as
the ``tilted grim-reaper'' and it can be parametrized as the  graph of the function $u(x,y)=\frac{\ln(\cos( y \cos \theta))}{\cos^2\theta}+x\tan \theta$ (here $\theta$ is a fixed  parameter that represents the angle of rotation and $x, y$ represent the coordinates of the two dimensional plane). Note that both the grim reaper and the tilted grim reaper are generated by a curve and are invariant in the other directions.\\

Another important translating solution is the ``bowl-soliton". This solution is a rotationally symmetric translating paraboloid that was obtained by Altschuler and Wu in \cite{AW}. It is the only noncollapsed and uniformly 2-convex translating solution (see \cite{Haslhofer}). More generally, rotationally symmetric translating solutions were classified by Clutterbuck, Schn\"urer and Schulze in \cite{CSS} and their stability was also studied,\\

Finally, we recall that by Wang \cite{Wang} every proper, convex translator in Euclidean space is either entire or lies in a slab region.

\subsection{Hyperbolic space: $\lambda_1 = \lambda_2 \neq 0$}

As in the Euclidean case, the study of the evolving curves in 2-dimensional hyperbolic space is better understood than the flow in higher dimensions. In \cite{NT2022} 
F. Nunes da Silva, K. Tenenblat classify soliton solutions to the curve shortening flow on two dimensional hyperbolic space. 
Unlike the classical case, many more solutions exist in the hyperbolic setting. Recall that in the classical case the grim reaper is the unique solution, up to rotations and translations, and in fact, in the hyperbolic case there exists a two parameter family of solutions to the equation.\\

Solitons in higher dimensional hyperbolic spaces were studied by G. Colombo, L.
Mari, M. Rigoli in 
\cite{CoMaRi}. Their work considers more generally  solitons in warped products, including  hyperbolic space as a special case of interest. They study solitons generated by  either Killing or (non-trivial) conformal vector fields $V$.  So, the problem reduces to solve the equation 
$$cV^{\bot} =\vec{ \bold{H}},$$
where $\vec{\bold{H}}$ is the mean curvature vector field.  Among their results they showed that
 the only entire graph $\psi : M \to \HH^{m+1} = \RR \times_{e^t} \RR^m$ in the $m+1$-dimensional hyperbolic space, which is a graphical soliton with respect to $V=e^t\partial_t$ with soliton constant $c < 0$ and contained in the half-space $\left[ \log\left(-\frac{2}{c}\right), + \infty\right)$ is the horosphere $\{t = \log\left(-\frac{2}{c}\right) \}$. More generally, they showed that
 there are no complete stable mean curvature flow solitons $\psi : M \to \HH^{m+1} = \RR\times_{e^t}\RR^m$ with respect to $V$ with soliton constant $c<0$ satisfying
	\begin{equation*}
	|\Phi| \in L^2(M) \quad \text{and} \quad \psi(M) \subseteq \left[ \log\left(-\frac{m}{c}\right), + \infty\right),
	\end{equation*}
	where $\Phi =  \bold{II} - \langle\;,\;\rangle_M \otimes \vec{\mathbf{H}}$ is the umbilicity tensor of $\psi$.
 \\

\subsection{The product space 
$\mathbb H^m \times  \mathbb R$: $\lambda_1 \neq  \lambda_2 = 0$.}

In this setting, Bueno studied rotationally symmetric translating solitons in \cite{Bueno}, proving the existence of
a solution that is a rotationally symmetric entire graph. This solution is analogous to the bowl soliton in Euclidean space and  is unique,  provided the asymptotic behavior at infinity coincides with a bowl lositon. Additionally, Bueno constructed a one parameter family of rotationally symmetric solutions with the topology of an annulus. These solutions are analogous to the wing-like solutions constructed by Clutterbuck, Schun\"urer and Schulze in \cite{CSS} for the Euclidean setting. In a subsequent paper  \cite{BO}, Bueno and Ortiz  extended  Bueno's results  to a related equation that prescribes mean curvature with a function of the normal vector, but these extensions do not introduce new solitons.\\

Analogues of the grim reaper and the tilted grim reaper  in this geometry have been classified by Lima and Pipoli in \cite{LP}. In that paper, the authors are concerned with a larger class of flows where the speed is given by symmetric polynomials on the principal curvatures. One of these polynomials agrees with mean curvature flow (the one that corresponds with $r=1$). In this broader context, the authors prove the existence translators and in the particular case of mean curvature flow they found  a one parameter family, called ``hyperbolic grim reapers''.

\subsection{ The Solvable group $Sol_3$: $\lambda_1 = - \lambda_2 = 1$.}

This case was studied by Pipoli in \cite{Pi1}. His  work  classifies solutions generated by a curve (and invariant in other directions); these could be considered as extensions of the grim reaper to this context. It is interesting to note that, in contrast with the Euclidean case (and Wang's result in \cite{Wang}), there are some solutions defined in slabs and others in half planes. Moreover, Pipoli proves the no-existence of non trivial solutions that are invariant in certain directions. Here  ``trivial''  refers to solutions that are minimal surfaces (classified in \cite{LM}). Note that rotational symmetry in Euclidean space implies that the choice of the invariant direction is irrelevant.

\section{Further solitons in $(S,g)$}\label{extensions}

In this section we will set up the necessary elements to state the equations of solitons in our family of $3$-dimensional solvmanifolds. 
To facilitate understanding for a broad audience, we include detailed computations for clarity.

\subsection{The geometry of submanifolds}\label{Sec:submfd}

Let $M$ denote an $n$-dimensional  embeded submanifold in $(\mathcal N,g)$ of codimension one. Let $TM$ denote its tangent space equipped with the induced metric. 
Thus,  around any point $p$ in $M$
$$T_p\mathcal N= T_pM \oplus \RR \nu,$$
where $\nu$ denotes a unit normal vector field to $T_pM$.\\
 
We now discuss the geometry on $M$ induced by the embedding, focusing on the definition of extrinsic geometric quantities. We remark that different sign conventions for these definitions exist in the literature and  we follow the one in \cite{Huisken-Polden}, which, in particular, differs from that in \cite{Pi1}.

We denote the induced metric (or first fundamental form) on $M$  by $\bar{g}$ and define it as
$$\overline{g}(X,Y)=g(X,Y) \hbox{ for } X,Y \in TM.$$
The second fundamental form is a  bilinear extrinsic geometric quantity that at any point $p \in M$ we denote by  $\alpha: T_pM\times T_p M\to\mathbb{R}$ and it is defined by
 $$\alpha(X,Y)=g(\nabla_X \nu,Y)=-g(\nu,\nabla_XY).$$
Here $\nabla$ denotes the Levi-Civita connection in $(\mathcal N,g)$ and $\nu$ a choice of unit normal vector field.

Note that since $g(\nu, \nu)=1$ we have that for every $X\in TM$,
$$ g(\nabla_X\nu, \nu)=0.$$
In consequence, the Weingarten map (also called shape operator) $A_\nu: TM\to TM $ is defined by
$$A_\nu (X)=\nabla_X\nu.$$
With this definition we have  $$\alpha(X,Y)=g(A_\nu X,Y)=- g(\nu,\nabla_XY).$$
In addition, for $X,Y$ vector fields on $M$ one has
$$\nabla_X Y = \overline{\nabla}_XY-\alpha(X,Y) \nu, $$
being $\overline{\nabla}$ the Levi-Civita connection on $M$.

 
 Since the Levi-Civita connection is torsion free we have
 $\nabla_UV-\nabla_VU=[U,V]$, for all $U,V\in T\mathcal N.$ Taking  vector fields $X,Y\in T\mathcal N$ that are tangent to $M$ and using the relationships above, we get
 $$[X,Y]=\nabla_XY-\nabla_YX=\overline{\nabla}_XY -\alpha(X,Y)\nu-\overline{\nabla}_YX + \alpha(Y,X)\nu.$$
 By looking at the orthogonal component one obtains that $\alpha$ is bilinear and symmetric. 
 This implies that $A_\nu$ is self-adjoint:
 $$g(A_{\nu}X,Y)=g(X, A_{\nu}Y).$$
The eigenvalues $\kappa_1, \ldots, \kappa_n$ of $A_\nu$ are known as the {\em  principal curvatures} and its trace is known as the {\em scalar mean curvature} (that we denote by $H$)
$$H=\kappa_1+\ldots+\kappa_n.$$
 We remark that if $\varphi_0:M\to \mathcal N$ is an embedding of $M$ into $\mathcal N$ we have that 
 $$\Delta_{M} \varphi_0=-H\nu,$$
where $\Delta_M$ is the Laplace-Beltrami operator.\\

 Given a basis $\{\tau_1,\ldots, \tau_n\}$ of $T_pM$, we introduce the following notation
\begin{align*}
\bar{g}_{ij}&= g(\tau_i, \tau_j),\\
(\bar{g}^{ij})_{ij}&=(\bar{g}_{ij})^{-1},\\
h_{ij} &=\alpha(\tau_i, \tau_j),\\
h^i_j&=\sum_{k=1}^n\bar{g}^{ik}h_{kj},\\
A_\nu(\tau_j)&= \sum_{i=1}^n h^i_j \tau_i.
\end{align*}

\subsubsection{The 2-dimensional case}
Following the notation above, when $M$ is a submanifold of dimension two, we have the following. 
 
 The matrix that represents the induced metric has coefficients
 $$E:=g(\tau_1,\tau_1),\qquad F:=g(\tau_1, \tau_2), \qquad G:=g(\tau_2, \tau_2),$$
 and they are ordered as follows
 \begin{equation}\label{metric}
\left( 
 \begin{matrix}	E & F \\
 	F & G
 	\end{matrix}
 \right).
\end{equation}
By making use of the formulas above, one has
$$  \left( \begin{matrix}
 h_{11} & h_{12} \\
h_{21}  & h_{22} 
\end{matrix} \right)   =-\left( \begin{matrix}
	g(\nabla_{\tau_1}\tau_1, \nu) & g(\nabla_{\tau_1}\tau_2, \nu)\\
	g(\nabla_{\tau_2}\tau_1, \nu) & g(\nabla_{\tau_2}\tau_2, \nu)
	\end{matrix} \right) = 
	\left( \begin{matrix}
		E & F\\
	F & G
	\end{matrix} \right)  
\left( \begin{matrix}
	h^1_{1}& h^1_{2}\\
	h^2_{1}& h^2_{2}
\end{matrix} \right).
$$
Here the matrix of the shape operator $A_{\nu}$ in the basis $\tau_1, \tau_2$  is the matrix
 \begin{equation}\label{shape}
	\left( 
	\begin{matrix}	h^1_{1} & h^1_{2} \\
		h^2_{1} & h^2_{2}
	\end{matrix}
	\right).
\end{equation}
Then we have the following formulas for the coefficients $h^i_{j}$:
$$\begin{array}{rcl}
	h^1_{1}& = & \frac1{EG-F^2}[Fg(\nabla_{\tau_2}\tau_1,\nu)-G g(\nabla_{\tau_1}\tau_1,\nu)],\\
	h^2_{1}& = & \frac1{EG-F^2}[F g(\nabla_{\tau_1}\tau_1,\nu)- Eg(\nabla_{\tau_2}\tau_1,\nu)],\\
	h^1_{2}& = & \frac1{EG-F^2}[-G g(\nabla_{\tau_1}\tau_2,\nu)+Fg(\nabla_{\tau_2}\tau_2,\nu)],\\
	h^2_{2}& = & \frac1{EG-F^2}[F g(\nabla_{\tau_1}\tau_2,\nu)-Eg(\nabla_{\tau_2}\tau_2,\nu)].
\end{array}
$$

	From these definitions, the trace of the shape operator $A_\nu$ is
\begin{align}\label{meanformula}
H  = & h^1_{1}+h^2_{2} \\
= &
\begin{array}{rcl}
		 \frac1{EG-F^2}[-G g(\nabla_{\tau_1}\tau_1,\nu)+Fg(\nabla_{\tau_2}\tau_1,\nu)+F g(\nabla_{\tau_1}\tau_2,\nu)-Eg(\nabla_{\tau_2}\tau_2,\nu)], \notag
		\end{array} \\
= &
\begin{array}{rcl}
		 \frac{-1}{EG-F^2}[-G h_{11} +Fh_{21}+F h_{12}-h_{22}], \notag
		\end{array}
\end{align}
	which corresponds to the mean curvature of the submanifold.

\subsection{Examples of submanifols in $3$-dimensional solvmanifolds}\label{ssec:Examples}

Here we will consider surfaces $M$ of $S$ and will calculate some of the operators described above. 
In Section \ref{ev of surfaces} we will use such estimates to obtain the mean curvature equation associated to these examples.

\subsubsection{$\widetilde{E}_1$- Invariant submanifolds}\label{ssec:ExamplesInvariant}


 Let $\gamma(s)=(0,y(s),z(s))$, $s \in I$ be a regular curve in $S$. We may assume that $\gamma$ is parametrized by arc-length, that is
\[
e^{-\lambda_2 z(s)} y'(s) = \cos \theta(s), \quad z'(s)=\sin\theta(s),
\]
where $\theta=\theta(s)$ is a smooth function.

Consider the surface   $M \subset S$  given by 
\[
\varphi(s,r)= L_{(r,0,0)}(\gamma(s))= (r,y(s),z(s)) \quad \textrm{ with } s\in I, \quad r\in \RR.
\]

It is easy to verify that $M$ is invariant by the translation maps $L_{(c,0,0)}$ which is the one-parameter subgroup for the Killing vector field $\widetilde{E}_1$. 
The tangent vectors  of $M$ are given by the vectors
\begin{align*}
\tau_1=\varphi_s(s,r)&=(0,y'(s),z'(s))=\cos(\theta) E_2 + \sin(\theta)E_3, \\
\tau_2=\varphi_r(s,r)&=(1,0,0)=e^{-\lambda_1 z} E_1,
\end{align*} where $\{E_1,E_2, E_3\}$ is the basis of left invariant vector fields described in Section \ref{sec:LieGroup}.

Then the unit normal associated to $M$ can be written as
 \[
 \nu= - \sin(\theta)E_2 + \cos(\theta)E_3.
 \] 
 The first fundamental form is given by the matrix 
 \begin{equation*}
\left( 
 \begin{matrix}	1 & 0 \\
 	0 &  e^{-2 \lambda_1 z}
 	\end{matrix}
 \right).
\end{equation*}

 We can also compute
 \begin{align*}
\nabla_{\tau_1}\tau_1&= (-\theta' \sin \theta-\lambda _2 \cos \theta \sin \theta)E_2+(\theta' \cos \theta+ \lambda_2 \cos^2 \theta) E_3\\&=(\theta' +\lambda_2 \cos \theta)(-\sin \theta E_2+\cos \theta E_3),
 \\  \nabla_{\tau_1}\tau_2 &= -\lambda_1 e^{-\lambda_1 z}\sin \theta E_1, 
  \\  \nabla_{\tau_2}\tau_1&=  -\lambda_1 e^{-\lambda_1 z}\sin \theta E_1, \\ \nabla_{\tau_2}\tau_2&=\lambda_1 e^{-2\lambda_1 z} E_3. 
\end{align*} 

 Then, to get the coefficients of the shape operator, we get
  \begin{align*}
h_{11}&
=-\theta'-\lambda_2\cos \theta,\\
  h_{12} &=0,
  \\  h_{21}&=0, \\ h_{22}&=  -\lambda_1 e^{-2\lambda_1 z} \cos \theta.
 \end{align*} 
In this case, the mean curvature (in Eq. \eqref{meanformula})
follows
\begin{equation}\label{mean1}
H=-\theta'(s)-(\lambda_1+\lambda_2)\cos \theta(s).
\end{equation}

\begin{rem} Notice that submanifolds of the form $\varphi(s,r)=(x(s), r, z(s))$ will be invariant by translations of the form $L_{(0,c,0)}$, and by following a reasoning as above, one gets the expression for the mean curvature $H$ as in Equation \eqref{mean1}.   For  $\gamma(s)=(x(s),0,z(s))$, $s \in I$ a regular curve in $S$ parametrized by arc-length, that is
\[
e^{-\lambda_2 z(s)} x'(s) = \cos \theta(s), \quad z'(s)=\sin\theta(s),
\]
where $\theta=\theta(s)$ is a smooth function.

\end{rem}


 \subsubsection{Graph submanifolds}\label{geometry of graphs}
 
 Now, we consider  submanifolds $M$ of $S$ that can be parametrized as $\varphi_0(y,z)=(f(y,z), y, z)$, where $f: \RR^2 \to \RR$ is a {smooth function}. 
 In this case, the tangent vectors are given by
 \begin{align*}
 \tau_1&= f_y  e^{-\lambda_1 z} E_1+ e^{-\lambda_2 z} E_2,\\
  \tau_2&= f_z e^{-\lambda_1 z} E_1 + E_3,
 \end{align*}
  where $f_y=\frac{\partial f}{\partial y}$, $f_z=\frac{\partial f}{\partial z}$. The unit normal associated to $M$ can be written as
 $$\nu= \frac{-E_1+f_y e^{(\lambda_2-\lambda_1)z}E_2+f_z  e^{-\lambda_1z}E_3}{ (1+f^2_y e^{2(\lambda_2-\lambda_1)z}+f^2_z  e^{-2\lambda_1z})^{\frac{1}{2}} }.$$
 The first fundamental form is given by
 \begin{equation*}
(\bar{g}_{ij})=
\left( 
 \begin{matrix}	 f^2_y  e^{-2\lambda_1 z} + e^{-2\lambda_2 z} & f_y f_z  e^{-2\lambda_1z} \\
 	 f_y f_z  e^{-2\lambda_1 z} &   f^2_z  e^{-2\lambda_1 z} + 1 
 	\end{matrix}
 \right).
\end{equation*}
 The inverse of this matrix is
  \begin{equation*} (\bar{g}^{ij})=
\frac{1}{f^2_y  e^{-2\lambda_1 z} +f^2_z  e^{-2(\lambda_1+\lambda_2) z} +  e^{-2\lambda_2 z}  }\left( 
 \begin{matrix}f^2_z  e^{-2\lambda_1 z} + 1 	& -f_y f_z  e^{-2\lambda_1 z} \\
 	 -f_y f_z  e^{-2\lambda_1 z} &   f^2_y  e^{-2\lambda_1 z} + e^{-2\lambda_2 z}  
 	\end{matrix}
 \right).
\end{equation*}
 We can also compute
 \begin{align*}
\nabla_{\tau_1}\tau_1&= f_{yy} e^{-\lambda_1 z} E_1+(f_y^2\lambda_1 e^{-2\lambda_1 z}+ \lambda_2 e^{-2\lambda_2 z} )E_3,
 \\ 
  \nabla_{\tau_1}\tau_2 &=  (f_{yz}-\lambda_1 f_y)e^{-\lambda_1 z}E_1-\lambda_2 e^{-\lambda_2 z} E_2+  f_z f_y\lambda_1 e^{-2\lambda_1 z}E_3, \\ 
   \nabla_{\tau_2}\tau_1&=(f_{yz}-\lambda_1 f_y)e^{-\lambda_1 z}E_1-\lambda_2 e^{-\lambda_2 z} E_2+  f_z f_y\lambda_1 e^{-2\lambda_1 z}E_3,  \\ 
   \nabla_{\tau_2}\tau_2&= (f_{zz}-2\lambda_1 f_z)e^{-\lambda_1 z} E_1+\lambda_1 f_z^2 e^{-2\lambda_1 z} E_3. 
\end{align*}

Then
\begin{align*}\label{ecugraph}
\sigma h_{11}&= -g(\nabla_{\tau_1}\tau_1, \nu) =[f_{yy}  -(f_y^2\lambda_1 e^{-2\lambda_1 z}+ \lambda_2 e^{-2\lambda_2 z} )f_z ]e^{-\lambda_1 z},\\
\sigma h_{12} &= -g(\nabla_{\tau_1}\tau_2, \nu)= [f_{yz} -f_y\lambda_1  +\lambda_2 f_y - f^2_z f_y\lambda_1 e^{-2\lambda_1 z}  ] e^{-\lambda_1 z},\\  
\sigma  h_{21}&= -g(\nabla_{\tau_2}\tau_1, \nu)= [f_{yz} -f_y\lambda_1  +\lambda_2 f_y -f^2_z f_y\lambda_1 e^{-2\lambda_1 z}  ] e^{-\lambda_1 z}, \\  
\sigma h_{22}&=  -g(\nabla_{\tau_2}\tau_2, \nu)= [f_{zz}-2\lambda_1 f_z -\lambda_1 f_z^3 e^{-2\lambda_1 z} ] e^{-\lambda_1 z}, 
\end{align*} 
here $\sigma=(1+f^2_y e^{2(\lambda_2-\lambda_1)z}+f^2_z  e^{-2\lambda_1z})^{\frac{1}{2}}$. 
 
\subsection{Setting up the equations for the construction of new solitons}\label{ev of surfaces}
  In this section we  consider solutions to Equation \eqref{eq:MCF} where $\{\varphi_t\}$ is a one-parameter family of isometries. 
  
  Let $V=\frac{\partial \varphi}{\partial t}$. If $V$ is a Killing vector field of the ambient space, then
  \eqref{eq:MCF} becomes
 \begin{equation}\label{eq-EvolutionKilling}
 g\left(V, \nu\right)=-H,
  \end{equation}
  and we say that $M$ is a translator in the direction of $V$. \\
  Below, we will describe the evolution of the two examples considered in Section \ref{ssec:Examples}, namely, 
solutions that have an additional invariance and solutions that are graphical in some direction, and discuss when they are translators.
In the future we would like to analize the existence and behavior of solutions of these examples.


\subsubsection{Evolution of $\widetilde{E}_1$- Invariant submanifolds}\label{ev of invariant surfaces}
 
Assume that $M$ is invariant under actions of $\widetilde{E}_1$. Then, M is parametrized as in  Section \ref{ssec:ExamplesInvariant}, and this parametrization satisfies the following system of differential equations:
\begin{equation*}
     \left\{
	       \begin{array}{l}
		 y'(s) = e^{\lambda_2 z(s)}  \cos \theta(s), \\
		 z'(s)=\sin\theta(s),  \\
		- \theta'(s)-(\lambda_1+\lambda_2)\cos \theta(s)=  H(s). 
	       \end{array}
	     \right.
   \end{equation*}\\

Let 
\[ 
V= \eta \widetilde{E}_1 + \beta \widetilde{E}_2 + \mu \widetilde{E}_3
\]
be a Killing vector field, where $\tilde E_i$ are the vector fields described in Section \ref{sec:LieGroup}. 
To see if $M$ is a translator in the direction of $V$ we can consider
$\eta=0$ since $\widetilde{E}_1$ is tangent to $M$. In this situation, 
\begin{align*}
\nu = & - \sin \theta E_2 + \cos \theta E_3\\ 
V= & \beta e^{- \lambda_2 z} E_2 + \mu (\lambda_1 x e^{-\lambda_1 z} E_1+ \lambda_2 y e^{-\lambda_2 z} E_2 + E_3).
\end{align*} 
It follows that $M$ is a translator in the direction of $V$ if and only if 
\[
-H = -\sin \theta \beta  e^{-\lambda_2 z} - \sin \theta \mu \lambda_2 y e^{-\lambda_2 z} + \cos \theta \mu.  
\] \\

Equivalently, 
\begin{align*}
\theta' &= -\sin \theta \beta  e^{-\lambda_2 z} - \sin \theta \mu \lambda_2 y e^{-\lambda_2 z} - \cos \theta (-\mu +\lambda_1+\lambda_2)\\
& =  - \beta z' e^{-\lambda_2 z} -  \mu \lambda_2 y z'e^{-\lambda_2 z} -y'  e^{-\lambda_2 z}  (-\mu +\lambda_1+\lambda_2)
\\& = \frac{d}{ds} \left[e^{-\lambda_2 z} \left(\frac{\beta}{\lambda_2}+\mu y\right)\right] -y'  e^{-\lambda_2 z}  (\lambda_1+\lambda_2) .
\end{align*}

Note that in the case considered in \cite{Pi1} $ \lambda_1+\lambda_2=0$, which simplifies the analysis of the ODE.

 \subsubsection{Evolution of graphs}\label{ev of graphs}

Assume that our surface is parametrized by $$\varphi_t(y,z)=(f(y,z)+t, y, z).$$ Take
 $V= \widetilde{E}_1 = e^{-\lambda_1 z} E_1$. 
 Then this surface is a soliton in the direction of $V$ if the following evolution equation holds  
  \begin{equation} 
  e^{-\lambda_1 z}=\sigma H,\label{graphev}
  \end{equation}
 where $\sigma= (1+f^2_y e^{2(\lambda_2-\lambda_1)z}+f^2_z  e^{-2\lambda_1z})^{\frac{1}{2}}$, $H=\sum_{i, j=1}^2 \bar{g}^{ij}h_{ji}$, 
and we have used \eqref{eq-EvolutionKilling} and the formulas obtained in Section \ref{geometry of graphs}.\\

 For simplicity, we now focus on solutions that are 1-dimensional, which are the analogous to the grim reaper in Euclidean space.  We first assume that $f$ only depends on $y$. We obtain
 \begin{equation}
   		\begin{array}{rcl}
   			   		\sigma h_{11} & = &  f_{yy}   e^{-\lambda_1 z},\\
					 \sigma  h_{12} & = & \sigma  h_{21}  =  (\lambda_2-\lambda_1)f_y e^{-\lambda_1 z},\\
   			   		\sigma h_{22} & = & 0.
   			   	\end{array}
   		   	\end{equation}
   		 The inverse of the first fundamental form is
  \begin{equation*} (\bar{g}^{ij})=
\frac{1}{f^2_y  e^{-2\lambda_1 z} +  e^{-2\lambda_2 z}  }\left( 
 \begin{matrix}     1 	& 0 \\
 	 0 &   f^2_y  e^{-2\lambda_1 z} + e^{-2\lambda_2 z}  
 	\end{matrix}
 \right) = \left( 
  \begin{matrix}  \tfrac{1}{f^2_y  e^{-2\lambda_1 z} +  e^{-2\lambda_2 z} } 	& 0 \\
 	 0 & 1
 	\end{matrix}
 \right).
\end{equation*}
	  Then equation \eqref{graphev} becomes 
\begin{equation*}\label{eq_f(y)}
 1= \frac{ f_{yy}}{f^2_y  e^{-2\lambda_1 z} +  e^{-2\lambda_2 z}  }.
 \end{equation*}

This equation cannot have a solution if $\lambda_1 \lambda_2 \ne 0$. Indeed, this can be justified as follows: Assume that there is a solutions $f$. 
 It is clear that constants cannot be solutions to the equation, hence
 $f_y\not \equiv 0$. Then we have for a suitable choice of $y$ that
$f^2_y  e^{-2\lambda_1 z} +  e^{-2\lambda_2 z} \to \infty$ either as $z\to \infty$ or as $z\to -\infty$ (depending on the signs of the  $\lambda_i$). This would imply that
  $1=0$,  which is a contradiction.

\bigskip
On the other hand, if we look at solutions that depend only on $z$ we obtain

 \begin{equation*} (\bar{g}^{ij})=
\left( 
 \begin{matrix} e^{2\lambda_2 z} 	& 0\\
 	0&   \frac{1}{f^2_z  e^{-2\lambda_1 z} + 1}  
 	\end{matrix}
 \right).
\end{equation*}

   	\begin{align*}
 \sigma h_{11}&=  - \lambda_2 e^{-(2 \lambda_2 +\lambda_1)z} f_z, \\
 \sigma h_{12} & =\sigma h_{21} =    0, \\ 
 \sigma h_{22}&=   [f_{zz}-2\lambda_1 f_z -\lambda_1 f_z^3 e^{-2\lambda_1 z} ] e^{-\lambda_1 z}. 
\end{align*}

  Then the \eqref{graphev} becomes
     	
$$1=\frac{  f_{zz}-2\lambda_1 f_z -\lambda_1 f_z^3 e^{-2\lambda_1 z}  }{ f^2_z  e^{-2\lambda_1 z} +1 }-  \lambda_2 f_z.
   $$

 In this case, standard ODE results provide local existence. We remark that in the Euclidean case there are one dimensional graphical translators in any direction, due to the invariance of the ambient metric with respect to rotations.\\

\noindent{\bf Acknowledgments.}
This project started at the ``Latin American and Caribbean Workshop on Mathematics and Gender'' held at Casa Matem\'atica Oaxaca (CMO) from May 15-20, 2022. There several groups were formed, and in particular, we were part of the Differential Geometry group. The groups were formed with the objective of creating collaboration among women mathematicians. \\
RA has been partially supported by CONICET, FONCYT and SeCyT-UNC.\\
RP is supported by CONACyT as an Investigadora por M\'exico and thanks IMUNAM Oaxaca for hosting her.\\
GO has been partially supported by FONCyT and SeCyT-UNR. \\
MS has been partially supported by Fondecyt Regular \#1190388.

\section{Declarations}

\noindent{\bf Conflict of interest.} No potential conflict of interest is reported by the authors.

\end{document}